\input amstex
\documentstyle{amsppt}

\magnification=1200
\pageheight{7.5in}
\expandafter\redefine\csname logo\string@\endcsname{}
\NoBlackBoxes

\topmatter
\title Topology of $3$-manifolds and a class of groups  
\endtitle
\author Sayed Khaled Roushon
\endauthor
\address School of Mathematics, Tata Institute, Homi Bhabha Road, Mumbai
400 005, India.
\endaddress
\email roushon\@math.tifr.res.in ,   
http://www.math.tifr.res.in/\~\ roushon/paper.html \endemail
\date November 24, 2002 \enddate
\thanks 2000 {\it Mathematics Subject Classification.}
Primary: 20F19, 57M99. Secondary: 20E99\endthanks
\abstract 
This paper grew out of an attempt to find a suitable finite sheeted
covering of an aspherical $3$-manifold so that the
cover either has infinite or trivial first homology group.
With this motivation we define a new class of groups. These groups are
in some sense eventually perfect. We prove results giving several classes
of examples of groups which do (not) belong to this class. Also we prove
some elementary results on these groups and
state two conjectures. A direct application of one of the conjectures to
the virtual Betti number conjecture of Thurston is mentioned.
\endabstract
\endtopmatter
\document
\baselineskip 14pt

\head 
0. Introduction
\endhead

The main motivation to this paper came from $3$-manifold topology  
while trying to find a suitable finite sheeted covering of an
aspherical $3$-manifold so that the cover has either infinite or trivial
first integral homology group. In \cite{R} it was proved that $M^3\times
{\Bbb D}^n$ is topologically rigid for $n>1$ whenever $H_1(M^3, {\Bbb Z})$ 
is infinite. Also the same result is true  when $H_1(M^3, {\Bbb Z})$ is
$0$. The remaining case is when $H_1(M^3, {\Bbb Z})$ is nontrivial
finite. There are induction techniques in surgery theory which can be
used to prove topological rigidity of a manifold if certain of finite
sheeted coverings of the manifold are also topologically rigid.
In the case of manifolds with nontrivial finite first integral homology
groups there is a natural finite sheeted cover namely the one which
corresponds to the commutator subgroup of the fundamental group.

So we start with a closed aspherical $3$-manifold $M$ with nontrivial
finite first integral homology group and consider the finite sheeted
covering $M_1$ of $M$ corresponding to the commutator subgroup. If
$H^1(M_1, {\Bbb Z})\neq 0$ or $H_1(M_1, {\Bbb Z})=0$ then we are done 
otherwise we again take the finite sheeted cover of $M_1$ corresponding 
to the commutator subgroup and continue. The group theoretic
conjecture (Conjecture 0.2) in this article implies that this process
stops in the sense that for some $i$ either $H^1(M_i, {\Bbb Z})\neq 0$ or 
$H_1(M_i, {\Bbb Z})=0$. 

Motivated by the above situation we define the following class of groups. 

\remark{Definition} An abstract group $G$ is called
{\it adorable} if $G^i/G^{i+1} = 1$ for some $i$, where
$G^i=[G^{i-1}, G^{i-1}]$, the commutator subgroup of $G^{i-1}$, and $G^0=G$. 
The smallest $i$ for which the above property is satisfied is called the 
{\it degree of adorability} of $G$. We denote it by $doa(G)$.\endremark

Obvious examples of adorable groups are finite groups, perfect groups,
simple groups and solvable groups. The first two classes are adorable
groups of degree $0$. The free products of perfect groups are adorable.
The abelian groups and symmetric groups on $n\geq 5$ letters are adorable
of degree $1$. Another class of adorable groups are $GL(R)=Lim_{n\to
\infty}GL_n(R)$. Here $R$ is any ring with unity and $GL_n(R)$ is the
multiplicative group of $n\times n$ invertible matrices. These are
adorable groups of degree $1$. This follows from the Whitehead lemma which
says that the commutator subgroup of $GL_n(R)$ is generated by the
elementary matrices and the group generated by the elementary matrices is a 
perfect group. Also $SL_n({\Bbb C})$, the multiplicative group of $n\times
n$ matrices with complex entries is a perfect group. In fact we will
prove that any connected Lie group is adorable as an abstract group. The
full braid groups on more than $4$ strings are adorable of degree $1$. 

We observe the following elementary facts in the next section:

\proclaim{Theorem 1.10} A group $G$ is adorable if and only if
there is a filtration $G_n< G_{n-1}<\cdots <G_1<G_0=G$ of
$G$ so that $G_i$ is normal in $G_{i-1}$, $G_{i-1}/G_i$ is abelian for
each $i$ and $G_n$ is a perfect group.\endproclaim

\proclaim{Theorem 1.14} Let $H$ be a normal subgroup of an adorable group
$G$. Then $H$ is adorable if one of the following conditions is 
satisfied: 
\roster
\item $G/H$ is solvable
\item for some $i$, $G^i/H^i$ is abelian
\item for some $i$, $G^i$ is simple
\item for some $i$, $G^i$ is perfect and the group $G^i/H^{i+1}$ does not
have any proper abelian normal subgroup
\endroster
\endproclaim

Also the braid group on more than $4$ strings are the examples to show
that an arbitrary finite index normal subgroup of an adorable group need
not be adorable.

\proclaim{Theorem 1.21} Every connected real or complex Lie group is
adorable as an abstract group.\endproclaim

Below we give some examples of non-adorable groups. Proofs of
non-adorability of some of these examples are easy. Proofs for the other
examples are given in the next sections.

Some examples of groups which are not adorable are non-abelian free
groups and fundamental groups of surfaces of genus greater than
$1$; for the intersection of a monotonically decreasing sequence of
characteristic subgroups of a non-abelian free group consists of the
trivial element only. The commutator subgroup of $SL_2({\Bbb Z})$ is
the nonabelian free group on 2 generators. Hence $SL_2({\Bbb Z})$ is not 
adorable. Also by Stallings' theorem fundamental groups of compact
$3$-manifolds which have finitely generated commutator subgroup with
infinite cyclic abelianization are also not adorable. It is known that 
most of these $3$-manifolds support hyperbolic metric by Thurston.
It is easy to show that the pure braid group is not adorable as there is
a surjection of any pure braid group of more than $2$ strings onto a
nonabelian free group.

The following results give some important classes of examples of
non-adorable groups.
 
\proclaim{Corollary 2.3} A torsion free Bieberbach groups is non-adorable
unless it is solvable.\endproclaim

\proclaim{Theorem 2.4} Let $G$ be a group satisfying the following
properties:
\roster
\item $H_1(G, {\Bbb Z})$ has rank $\geq 3$ 
\item $H_2(G^j, {\Bbb Z})=0$ for $j\geq 0$
\endroster

Then $G$ is not adorable. Moreover, $G^j/G^{j+1}$ has rank $\geq 3$ for
each $j\geq 1$.\endproclaim

The proposition below is a consequence of the above Theorem.

\proclaim{Proposition 2.7} A knot group is adorable if and
only if it has trivial Alexander polynomial.\endproclaim

In fact in this case the commutator subgroup of the knot group is perfect.
All other knot groups are not adorable. On the other hand any knot
complement supports a complete nonpositively curved Riemannian metric
({\cite{L}). 

After seeing an earlier version of this paper (\cite{R2}) Tim Cochran
informed me that the Proposition 2.7 was also observed by him in
[Corollary 4.8, \cite{C}]. 

Note that all the torsion free examples of non-adorable groups we
mentioned above act freely and properly discontinuously (except the braid
groups case, which is still an open question) on a simply connected
complete nonpositively curved Riemannian manifold. Also we recall
that a solvable subgroup of the fundamental group of a nonpositively
curved manifold is virtually abelian (\cite{Y}). There are
generalization of these results to the case of locally $CAT(0)$ spaces
(\cite{BH}). Considering these facts we pose the following conjecture.

\proclaim{Conjecture 0.1} Fundamental group of generic class
of complete nonpositively curved Riemannian manifolds or more generally of
generic class of locally $CAT(0)$ metric spaces are not
adorable.\endproclaim

One can even ask the same question for hyperbolic groups.

Now we state the conjecture we referred before.   

\proclaim{Conjecture 0.2} Let $G$ be a finitely presented torsion free 
group such that $G^i/G^{i+1}$ is a finite group for all $i$.
Then $G$ is adorable.\endproclaim 

Note that $G^i/G^{i+1}$ is finite for each $i$ if and only if $G/G^i$
is finite for each $i$. Thus the above conjecture says that a non-adorable
finitely presented torsion free group has an infinite solvable quotient.
Compare this observation with Proposition 2.1. 

There is another consequence of this conjecture. That is, 
proving this conjecture for the particular case when the group $G$ is the 
fundamental group of an aspherical $3$-manifold will imply that the
virtual Betti number conjecture of Thurston is true if a modified (half)
version of it is true. We mention it below:

\proclaim{Modified virtual Betti number conjecture} Let $M$ be a
closed 
aspherical $3$-manifold such that $H_1(M, {\Bbb Z})=0$. Then there is a 
finite sheeted covering $\tilde M$ of $M$ with $H_1(\tilde M, {\Bbb Z})$ 
infinite.\endproclaim

It is easy to see that the Conjecture 0.2 and the Modified virtual Betti
number conjecture together implies the virtual Betti number conjecture.

\proclaim{Virtual Betti number conjecture} Any closed aspherical
$3$-manifold has a finite sheeted covering with infinite first homology
group.\endproclaim

\remark{Acknowledgment} Part of this work was presented in the conference
on Algebraic and Geometric Topology, January 01-04, 2002, Delhi
University, India and in the
Satellite conference on Geometric Topology of the ICM2002, August 12-16, 
2002, Shaanxi Normal
University, Xi'an, China. The author would like to thank the
organizing committees for the invitation to participate and lecture in the
conferences.\endremark

\head 
1. Elementary facts about Adorable groups
\endhead

In this section we prove some basic results on adorable groups. 

Recall that a group is called {\it perfect} if the commutator subgroup of
the group is the whole group.

\proclaim{Proposition 1.1} Let $f:G\to H$ be a surjective homomorphism
with $G$ adorable. Then $H$ is also adorable and $doa(H)\leq
doa(G)$.\endproclaim

\demo{Proof} $f$ induces surjective homomorphism $G^i\to H^i$ for each
$i$. The proof follows from the definitions of adorable groups and
its' $doa$.\qed\enddemo

\remark{Example 1.2} The Artin pure braid group on more than $2$ strings
is not adorable, for it has a quotient a non-abelian free group. 
In fact the full braid group on $n$-strings is not adorable for
$n\leq 4$ and adorable of degree $1$ otherwise. (see \cite{GL}).\endremark 

\proclaim{Proposition 1.3} The product $G\times H$ of two groups are
adorable if and only if both the groups $G$ and $H$ are adorable. Also if
$G\times H$ is adorable then $doa(G\times H)=\text{max}\ \{doa(G),
doa(H)\}$.\endproclaim

\demo{Proof} Note that $(G\times H)^i=G^i\times H^i$. If $G$ and $H$ are both 
adorable then it clearly follows that so is their product and also it follows 
that $doa(G\times H)=\text{max}\ \{doa(G), doa(H)\}$. The `only if' part 
follows from Proposition 1.1.\qed\enddemo

\remark{Remark 1.4} Note here that free product of two adorable groups 
need not be adorable; for example the non-abelian free group on two 
generators is not adorable, but the infinite cyclic group is.
But it is plausible that the free product of two nonsolvable adorable
group is adorable.\endremark

\proclaim{Proposition 1.5} Let $H$ be a normal subgroup of a group $G$
with quotient $F$ such that both $H$ and $F$ are perfect, then $G$ is also
perfect.\endproclaim

\demo{Proof} We have $$G/H=F=F^1=(G/H)^1=G^1/{G^1\cap H}=G^1/{G^1\cap
H^1}=G^1/H^1=G^1/H$$ This proves the Proposition.\qed\enddemo

\proclaim{Lemma 1.6} Let $G$ be an adorable group and $H$ is a normal 
subgroup of $G$. Assume that for some $i_0$, $G^{i_0}$ is simple. Then 
$H$ is also adorable and $doa(H)\leq doa(G)$.\endproclaim

\demo{Proof} Note that $H^{i_0}$ is a normal subgroup of $G^{i_0}$ and 
hence either $H^{i_0}=1$ or $H^{i_0}=G^{i_0}$. In any case $H$ is 
adorable and $doa(H)\leq doa(G)$.\qed\enddemo

\remark{Remark 1.7} In the above lemma instead of assuming the strong 
hypothesis that $G^{i_0}$ is simple we can assume only that $G^{i_0}$ 
is perfect and $G^{i_0}/H^{i_0+1}$ does not have any proper normal abelian
subgroup. With this weaker hypothesis the proof follows from the fact that
the kernel of the surjective homomorphism $G^{i_0}/H^{i_0+1}\to
G^{i_0}/H^{i_0}$ is either trivial or $G^{i_0}=H^{i_0}$. In either case 
it follows that $H$ is adorable.\endremark 

\proclaim{Lemma 1.8} Let $H$ be a normal subgroup of an adorable group 
$G$ such that $G^i/H^i$ is abelian for some $i$. Then $H$ is also
adorable.\endproclaim

\demo{Proof} There is an $i_0>i$ so that $G^{i_0+1}=G^{i_0}$. Also as 
$G^i/H^i$ is abelian we get $G^{i+1}\subset H^i$. Now we have
$$H^{i_0}\subset 
G^{i_0}=G^{i_0+2}=G^{i+1+(i_0-i+1)}\subset H^{i_0+1}$$ Also
$H^{i_0+1}\subset H^{i_0}$. 
Hence $H^{i_0+1}=H^{i_0}$. Therefore $H$ is adorable.\qed\enddemo

\proclaim{Lemma 1.9} Let $H$ be a normal subgroup of a group $G$ such
that 
$G^i/H^i$ is abelian for some $i$. Then $G$ is adorable if and only if $H$
is adorable.\endproclaim 

\demo{Proof} One way of the proof follows from the previous Lemma. So
assume that $H$ is adorable. Choose $i_0>i$ so that $H^{i_0+1}=H^{i_0}$. 
From the hypothesis $G^{i+1}\subset H^i$. Now
$$G^{i_0+1}=G^{i+1+i_0-i}\subset H^{i_0}=H^{i_0+2}\subset G^{i_0+2}$$ 
Also we have $G^{i_0+2}\subset G^{i_0+1}$. Thus $G^{i_0+2}=G^{i_0+1}$ 
and hence $G$ is adorable.\qed\enddemo

\proclaim{Theorem 1.10} A group $G$ is adorable if and only if    
there is a filtration $G_n< G_{n-1}<\cdots <G_1<G_0=G$ of
$G$ so that $G_i$ is normal in $G_{i-1}$, $G_{i-1}/G_i$ is abelian for
each $i$ and $G_n$ is a perfect group.\endproclaim 

\demo{Proof} We use Proposition 1.11 below and induction on $n$ to prove
the `if' part of the Theorem. So assume that there is a filtration of $G$
as in the hypothesis. Then there is an exact sequence $$1\to
G_n\to G_{n-1}\to G_{n-1}/G_n\to 1$$ such that $G_{n-1}/G_n$ is abelian
and $G_n$ is perfect and hence adorable. By Proposition 1.11 $G_{n-1}$ is
adorable. By induction $G_0=G$ is adorable. The `only if' part 
of the Theorem follows from the definition of adorable
groups.\qed\enddemo

\proclaim{Proposition 1.11} Let $H$ be a normal subgroup of a group 
$G$ such that $G/H$ is solvable. Then $H$ is adorable if and only if so is
$G$.\endproclaim

\demo{Proof} Before we start with the proof, we note down some generality: 
Suppose $G$ has a filtration as in the hypothesis of Theorem 1.10.  
Since $G_{i-1}/G_i$ is abelian for each $i$, we have ${G'}_{i-1}\subset
G_i$. Replacing $i$ by $i+1$ we get $G_i'\subset G_{i+1}$. 
Consequently, $G_0^i=G^i=\{G'\}^{i-1}\subset {G_1}^{i-1}\subset {\{G_1}'\}^{i-2}
\subset {G_2}^{i-2}\subset\cdots \subset {G'}_{i-1}\subset G_i$. Thus we get 
$G^n\subset G_n$. 

Denote $G/H$ by $F$. As $F$ is solvable we have $1\subset F^k\subset\cdots 
\subset F^1\subset F^0=F$ where $F^k$ is abelian. Let $\pi:G\to G/H$ be
the quotient map. We have the following sequence of normal subgroups of
$G$: $$\cdots\subset H^n\subset H^{n-1}\cdots\subset H^1\subset H\subset
\pi^{-1}(F^k)\cdots\subset \pi^{-1}(F^0)=G$$

Note that this sequence of normal subgroup satisfies the same properties 
as those of the filtration $G_i$ of $G$ above. Hence $G^{k+i}\subset
H^{i-1}$. Now if $G$ is adorable then for some $i$, $G^{k+i}$ is
perfect. We have $$H^{k+i}\subset G^{k+i}=G^{k+k+i+2}\subset
H^{k+i+1}$$ But we already have $H^{k+i+1}\subset H^{k+i}$. That is
$H^{k+i}$ is perfect, hence $H$ is adorable. Conversely if $H$ is adorable
then for some $i$, $H^i$ is perfect. By Theorem 1.10 it follows that
$G$ is also adorable.\qed\enddemo 

\proclaim{Corollary 1.12} Let $G$ be a torsion free infinite group and
$F$ be a finite quotient of $G$ with kernel $H$ such that $H$ is free
abelian and also central in $G$. Then $G$ is adorable.\endproclaim

\demo{Proof} Recall that equivalence classes of extensions of $F$ by 
$H$ are in one to one correspondence with $H^2(F, H)$ which is isomorphic
to $Hom(F, ({\Bbb R}/{\Bbb Z})^n)$ where $n$ is the rank of $H$ (see
exercise 3, page 95, in \cite{Br}). 
If $F$ is perfect then $Hom(F, ({\Bbb R}/{\Bbb Z})^n)=0$ and hence the  
extensions $1\to H\to G\to F\to 1$ splits. But by hypothesis $G$ is
torsion free. Hence $F$ is not perfect. By a similar argument it can be
shown that $F^i$ is perfect for no $i$ unless it is the trivial group.
Since $F$ is finite this proves that $F$ is solvable and hence $G$ is
adorable, in fact solvable.\qed\enddemo

We sum up the above Lemmas and Propositions in the following Theorem.

\proclaim{Theorem 1.13} Let $H$ be a normal subgroup of an adorable group
$G$. Then $H$ is adorable if one of the following conditions is 
satisfied:
\roster
\item $G/H$ is solvable
\item for some $i$, $G^i/H^i$ is abelian
\item for some $i$, $G^i$ is simple
\item for some $i$, $G^i$ is perfect and the group $G^i/H^{i+1}$ does not
have any proper abelian normal subgroup
\endroster
\endproclaim

\remark{Remark 1.14} It is known that any countable group is a subgroup of
a countable simple group (see theorem 3.4, chapter IV in \cite{LS}). Also
we mentioned before that even finite index normal subgroup of an adorable
group need not be adorable. So the above theorem is best possible in this
regard.
\endremark 

In the next section we give some more examples of virtually adorable
groups which are not adorable. 

The following is an analogue of a theorem of Hirsch for poly-cyclic 
groups. 

\proclaim{Theorem 1.15} The following are equivalent:

\roster
\item $G$ is a group which admits a filtration $G=G_0>G_1>\cdots >G_n$
with the property that each $G_{i+1}$ is normal in $G_i$ with quotient
$G_i/G_{i+1}$ cyclic and $G_n$ is a perfect group which satisfies the
maximal condition for subgroups.

\item $G$ is adorable and it satisfies the maximal condition for 
subgroups, i.e., for any sequence $H_1<H_2<\cdots$ of subgroups of 
$G$ there is an $i$ such that $H_i=H_{i+1}=\cdots$.

\endroster\endproclaim

\demo{Proof} The proof is on the same line as Hirsch's theorem. The main
lemma is the following:

\proclaim{Lemma A} Let $H_1$ and $H_2$ be two subgroup of a group $G$ and 
$H_1\subset H_2$. Let $H$ be a normal subgroup of $G$ with the 
property that $H\cap H_1=H\cap H_2$ and the subgroup generated by $H$ and
$H_1$ is equal to the subgroup generated by $H$ and $H_2$. Then
$H_1=H_2$.\endproclaim  

{\bf (1) implies (2):} By Theorem 1.10 it follows that $(1)$ implies 
that $G$ is adorable. Now we check the maximal condition by induction 
on $n$. As $G_n$ already satisfy maximal condition we only need to 
check that $G_{n-1}$ also satisfy maximal condition which follows from 
the following Lemma and by noting that $G_{n-1}/G_n$ is cyclic: 

\proclaim{Lemma B} Let $H$ be a normal subgroup of a group $G$ such that 
both $H$ and $G/H$ satisfy the maximal condition then $G$ also satisfies
the maximal condition.\endproclaim

\demo{Proof} Let $K_1<K_2<\cdots$ be an increasing sequence of subgroups
of $G$. Consider the two sequences of subgroups $H\cap K_1<H\cap
K_2<\cdots$ and $\{H, K_1\}<\{H, K_2\}<\cdots$. Here $\{A, B\}$ denotes 
the subgroup generated by the subgroups $A$ and $B$. As $H$ and $G/H$ both
satisfy the maximal condition there are integers $k$ and $l$ so that 
$H\cap K_k=H\cap K_{k+1}=\cdots$ and $\{H, K_l\}=\{H, K_{l+1}\}=\cdots$. 
Assume $k\geq l$. Then by Lemma A $K_k=K_{k+1}=\cdots$.\qed\enddemo

{\bf (2) implies (1):} As $G$ is adorable it has a filtration
$G=G_0>G_1>\cdots >G_n$ with $G_n$ perfect and each quotient abelian. 
Also $G_n$ satisfies maximal condition as it is a subgroup of $G$ and $G$
satisfies maximal condition. Since $G$ satisfies maximal condition 
each quotient $G_i/G_{i+1}$ is finitely generated. Now a filtration as in 
$(1)$ can easily be constructed. 

This proves the theorem.\qed\enddemo

\head
2. Some examples of (non-)adorable groups
\endhead

This section gives a large and important class of examples of
non-adorable groups.

\proclaim{Proposition 2.1} Let $M^3$ be a compact $3$-manifold with the 
property that there is an exact sequence of groups $1\to H\to \pi_1(M)\to
F\to 1$ such that $H$ is finitely generated nonabelian but not the
fundamental group of the Klein bottle and $F$ is an
infinite solvable group. Then $\pi_1(M)$ is not adorable.\endproclaim

\demo{Proof} By Theorem 11.1 in \cite{He} it follows that $H$ is the
fundamental group of a compact surface. Also as $H$ is not the Klein
bottle group, it is not adorable. The Corollary now follows from 
Proposition 1.11.\qed\enddemo

\proclaim{Proposition 2.2} Let $G$ be a torsion free group and $H$ a free
(abelian or non-abelian) normal subgroup of $G$ with quotient $F$ a
non-trivial finite perfect group. Then $G$ is not adorable.\endproclaim

\demo{Proof} If $H$ is non-abelian then by Stallings' Theorem $G$
itself is free and hence not adorable. So assume $H$ is free abelian. 
Since $F$ is a perfect group, the restriction of the quotient
map $G\to F$ to $G^i$ is again surjective for each $i$ with $H\cap G^i$ as
kernel. And since $G$ is infinite and torsion free, $H\cap G^i$ is
non-trivial free abelian for all $i$. This shows that each $G^i$ is
again a Bieberbach group. It is known that if $H^1(G^i, {\Bbb Z})=0$ then 
$G^i$ is centerless and centerless Bieberbach groups are meta-abelian and
hence solvable (\cite{HS}). But since each $G^i$ surjects 
onto a non-trivial perfect group it cannot be solvable. Hence 
$H^1(G^i, {\Bbb Z})\neq 0$ for each $i$. This proves the
Proposition.\qed\enddemo

The conclusion of the above Proposition remains valid if we assume that
$F$ is non-solvable adorable.

By Bieberbach theorem (\cite{Ch}) we have the following corollary.

\proclaim{Corollary 2.3} The fundamental group of a closed flat
Riemannian manifold is non-adorable unless it is solvable.\endproclaim

The following Theorem gives some more examples of non-adorable groups.

\proclaim{Theorem 2.4} Let $G$ be a group satisfying the following
properties:
\roster
\item $H_1(G, {\Bbb Z})$ has rank $\geq 3$ 
\item $H_2(G^j, {\Bbb Z})=0$ for $j\geq 0$
\endroster

Then $G$ is not adorable. Moreover, $G^j/G^{j+1}$ has rank $\geq 3$ for
each $j\geq 1$.\endproclaim

\demo{Proof} Consider the short exact sequence. $$1\to G^1\to G\to
G/G^1\to 1$$

We use the Hochschild-Serre spectral sequence (\cite{Br}, page 171) of the
above exact sequence. The $E^2$-term of the spectral sequence is  
$E_{pq}^2=H_p(G/G^1, H_q(G^1, {\Bbb Z}))$. Here ${\Bbb Z}$ is
considered as a trivial $G$-module. This spectral sequence gives rise to
the following five term exact sequence. $$H_2(G, {\Bbb Z})\to
E^2_{20}\to E^2_{01}\to H_1(G, {\Bbb Z})\to E^2_{10}\to 0$$ 

Using $(2)$ we get $$0\to H_2(G/G^1, H_0(G^1, {\Bbb Z}))\to H_0(G/G^1,
H_1(G^1, {\Bbb Z}))\to H_1(G, {\Bbb Z})\to$$$$\to H_1(G/G^1, H_0(G^1,
{\Bbb Z}))\to 0$$ 

As $\Bbb Z$ is a trivial $G$-module we get $$0\to H_2(G/G^1, {\Bbb Z})\to
H_0(G/G^1, H_1(G^1, {\Bbb Z}))\to H_1(G, {\Bbb Z})\to H_1(G/G^1, {\Bbb 
Z})\to 0$$

Note that the homomorphism between the last two non-zero terms in the
above exact sequence is an isomorphism. Also the second non-zero term from
left is isomorphic to the co-invariant $H_1(G^1, {\Bbb Z})_{G/G^1}$ and
hence we have the following $$H_2(G/G^1, {\Bbb Z})\simeq H_1(G^1, {\Bbb
Z})_{G/G^1}$$ 

Since $G/G^1$ has rank $\geq 3$ we get that $H_2(G/G^1, {\Bbb Z})$ has
rank greater or equal to $^3C_2=3$. This follows from the following
lemma. 

\proclaim{Lemma 2.5} Let $A$ be an abelian group. Then the rank
of $H_2(A, {\Bbb Z})$ is $^{rk A}C_2$ if $rk A$ is finite otherwise it
is infinity.\endproclaim

\demo{Proof} If $A$ is finitely generated then from the formula
$H_2(A, {\Bbb Z})\simeq \bigwedge^2 A$ it follows that rank of $H_2(A,
{\Bbb Z})$ is $^{rk A}C_2$. In the case $A$ is countable and infinitely
generated then there are finitely generated subgroups $A_n$ of $A$  such
that $A$ is the direct limit of $A_n$. Now as homology of group
commutes with direct limit the proof follows using the previous
case. Similar argument applies when $A$ is uncountable.\qed\enddemo  

To complete the proof of the theorem note that there is a surjective
homomorphism $H_1(G^1, {\Bbb Z})\to H_1(G^1, {\Bbb Z})_{G/G^1}$. Thus we
have proved that $H_1(G^1, {\Bbb Z})$ also has rank $\geq 3$. Finally
replacing $G$ by $G^n$ and $G^1$ by $G^{n+1}$ and using induction on $n$
the proof is completed.\qed\enddemo

There are two important consequences of Theorem 2.4. At first we recall
some definition from \cite{St}. 

Let $R$ be a non-trivial commutative ring with unity. The class $E(R)$
consists of groups $G$ for which the trivial $G$-module $R$ has a
$RG$-projective resolution $$\cdots\to P_2\to P_1\to P_0\to R\to 0$$ 
such that the map {\bf$1$}$_R\otimes \partial_2:R\otimes_{RG}P_2\to
R\otimes_{RG}P_1$ is injective. Note that if a group belongs to  $E(R)$
then $H_2(G, R)=0$. Also this condition is sufficient to belong to
$E(R)$ for groups of cohomological dimension less or equal to $2$. By
definition $G$ lies in $E$ if it belongs to $E(R)$ for all $R$. A 
characterization of $E$-groups is that a group $G$ is an $E$-group if and
only if $G$ belongs to $E({\Bbb Z})$ and $G/G^1$ is torsion free (lemma
2.3 in \cite{St}).

\proclaim{Corollary 2.6} Let $G$ be an $E$-group and rank of $H_1(G, {\Bbb
Z})$ is $\geq 2$. Then $G$ is not adorable.\endproclaim

\demo{Proof} By theorem A in \cite{St} it follows that $G$ satisfies
condition $(2)$ of Theorem 2.4. From the proof of Theorem 2.4 we get 
that $H_1(G^2, {\Bbb Z})$ has rank $\geq 1$ and hence in particular $G^2$
is not perfect. On the other hand an $E$-groups has derived length $0,1,2$
or infinity (remark after theorem A in \cite{St}). Thus $G$ is not 
adorable.\qed\enddemo 

In the following Proposition we give an application of the above Theorem
for knot groups.

\proclaim{Proposition 2.7} Let $H=\pi_1({\Bbb S}^3-k)$, where $k$ is a 
nontrivial knot in the $3$-sphere with non-trivial Alexander polynomial.
Then $H$ is not adorable. Moreover if rank of $H^1/H^2$ is
greater or equal to $3$ then the same is true for
$H^j/H^{j+1}$ for all $j\geq 2$.\endproclaim

In fact a stronger version of the Proposition follows, namely by \cite{St}  
the successive quotients of the derived series of $G$ are
torsion free. Thus we get that the successive quotients of the derived
series are nontrivial and torsion free. 

\demo{Proof of Proposition 2.7} At first recall that condition $(2)$ of
Theorem 2.4 follows from theorem A in \cite{St}. On the other hand the
commutator subgroup of a knot group is perfect if and only if the knot has
trivial Alexander polynomial. So assume that $H^1$ is not perfect. If
$H^1$ is finitely generated then in fact it is non-abelian free and hence
$H$ is not adorable. If rank of $H^1/H^2$ is $\geq 3$ then the proof
follows from the above Theorem. So assume that rank of $H^1/H^2$ is $\leq
2$. 

Recall that the rank of the abelian group $H^1/H^2$ is equal to the degree
of the Alexander polynomial of the knot (see theorem 1.1 in \cite{Cr}).
Thus if rank of $H^1/H^2$ is $1$ then the Alexander polynomial has degree
$1$ which is impossible as the Alexander polynomial of a knot always has
even degree. Next if rank of $H^1/H^2$ is $2$ then $H$ is not adorable by
Corollary 2.6 and noting that knot groups are $E$-groups.\qed\enddemo 

\proclaim{Definition 2.7} {\rm A Lie group is called {\it adorable} if it
is adorable as an abstract group.}\endproclaim

\proclaim{Theorem 2.8} Every connected (real or complex) Lie group is
adorable.\endproclaim

\demo{Proof} Let $G$ be a Lie group and consider its derived 
series: $$\cdots\subset G^n\subset G^{n-1}\cdots \subset G^1\subset
G^0=G$$ 

Note that each $G^i$ is a normal subgroup of $G$. Define
$G_i=\overline{G^i}$. Then we have a sequence of normal subgroups: 
$$\cdots\subset G_n\subset G_{n-1}\cdots \subset G_1\subset G_0=G$$ 
so that $G_i$ is a closed Lie subgroup of $G$ and $G_i/G_{i+1}$ is
abelian for each $i$. Suppose for some $i$, dim $G_i=0$,  
i.e., $G_i$ is a closed discrete normal subgroup of $G$. We claim 
$G_i$ is abelian. For, fix $g_i\in G_i$ and consider the continuous map
$G\to G_i$ given by $g\mapsto gg_ig^{-1}$. As $G$ is connected 
and $G_i$ is discrete image of this map is the singleton $\{g_i\}$. That
is $g_i$ commutes with all $g\in G$ and hence $G_i$ is abelian. 

As $G^i\subset G_i$, $G^i$ is also abelian. Thus $G$ is solvable and hence 
adorable.

Next assume no $G_i$ is discrete. Then as $G$ is finite dimensional and
$G_i$'s are Lie subgroup of $G$  
there is an $i_0$ so that $G_j=G_{j+1}$ for all $j\geq i_0$ and
dim $G_{i_0}\geq 1$. We need the following Lemma to complete
the proof of the Theorem.

\proclaim{Lemma 2.9} Let $G$ be a (real or complex) Lie group
such that $\overline {G^1}=G$. Then $G^2=G^1$, that is $G^1$ is a perfect
group.\endproclaim

\demo{Proof} The proof follows from Theorem XII.3.1 and Theorem XVI.2.1 of 
\cite{Ho}.\qed\enddemo

We have $G^{i_0}\subset G_{i_0}$ and hence $$G_{i_0}=G_{i_0+1}
=\overline{G^{i_0+1}}\subset \overline{G^1_{i_0}}\subset
\overline{G_{i_0}}=G_{i_0}$$ This implies $\overline{G^1_{i_0}}=G_{i_0}$. 
Now from the above Lemma we get $G_{i_0}$ is adorable. Thus $G_{i_0}$ is a
normal adorable subgroup of $G_{i_0-1}$ with quotient $G_{i_0-1}/G_{i_0}$ 
abelian and hence by Proposition 1.11 $G_{i_0-1}$ is also adorable. By
induction it follows that $G$ is adorable.\qed\enddemo

\enddocument